# ON THE SUBGROUPS OF THE EXTENDED HECKE, HECKE AND PICARD GROUPS

Ma. Louise Antonette N. De Las Peñas
Mathematics Department, Ateneo de Manila University, Loyola Heights, Quezon City, 1108 Philippines
e-mail: mlp@math.admu.edu.ph

Ma. Carlota B. Decena
Department of Mathematics and Physics, College of Science, University of Santo Tomas, Espana, Manila, Philippines
e-mail: cbdecena@mnl.ust.edu.ph

Glenn R. Laigo
Mathematics Department, Ateneo de Manila University, Loyola Heights, Quezon City, 1108 Philippines
e-mail: glaigo@ateneo.edu

**Abstract.** In this work, we derive the low index subgroups of the extended Hecke, Hecke and the Picard groups using tools in color symmetry theory. We also present the low index subgroups of the modular group $23\infty$.

## 1 Introduction

In [8] Erich Hecke introduced an infinite class of discrete subgroups $H(\lambda_p)$ of $PSL(2, R)$ consisting of linear fractional transformations preserving the upper half plane. The Hecke group $H(\lambda_p)$ is the group generated by $X(z) = -1/z$ and $Y(z) = -1/(z + \lambda_p)$ where $\lambda_p = 2\cos(\pi/p)$, $p \in N$, $p \geq 3$ having presentation $<X, Y : X^2 = Y^p = e>$. Adding the reflection $R(z) = 1/\bar{z}$ to the generators $X$, $Y$ of $H(\lambda_p)$ gives the extended Hecke group $\overline{H}(\lambda_p)$ with presentation $<R, X, Y : R^2 = X^2 = Y^p = (RX)^2 = (YR)^2 = e>$. The extended Hecke group can also be generated by $P$, $Q$ and $R$ where $X = RP$ and $Y = QR$ with defining relations $P^2 = Q^2 = R^2 = (RP)^2 = (QR)^p = e$.

Geometrically, $P$, $Q$ and $R$ denote reflections along the sides of a triangle $\Delta$ with interior angles $\pi/2$, $\pi/p$ and zero. In this sense, the extended Hecke group is considered a special type of triangle group denoted by $*2p\infty$. The Hecke group $H(\lambda_p)$ is a subgroup of index 2 in $\overline{H}(\lambda_p)$ and is also referred to as $2p\infty$. A well known Hecke group is the group $H(\lambda_3)$ or $23\infty$. It is more known as the modular group $PSL(2, Z)$, consisting of linear fractional transformations with integral coefficients.

Hecke groups [6, 9, 18, 21], extended Hecke groups [12, 17, 19] and the Picard group [6, 7, 10, 16, 20, 21] have been studied extensively for many aspects in the literature. The modular group, as well as the Picard group has especially been of great interest in many fields of mathematics, for example, number theory, automorphic function theory and group theory. The modular group, in particular, plays a crucial role in the uniformization of algebraic curves and in the construction of subgroups of larger matrix groups (e.g, $SL(n, Z)$ and $SL(n, R)$). The Picard group contains the modular group as a subgroup.

In this work, the main focus would be to determine the low index subgroups of the extended Hecke, Hecke and Picard groups using a geometric approach based on tools in color symmetry theory. For an extended Hecke group $*2p\infty$ or a Hecke group $2p\infty$, deriving its subgroups entails treating the group as a group of symmetries of a two dimensional tiling by a hyperbolic triangle $\Delta$ with interior angles $\pi/2$, $\pi/p$ and zero.

A similar setting is assumed to derive the low index subgroups of the Picard group or the Gaussian Modular group $PSL(2, Z[i])$, the three dimensional counterpart of the modular group, consisting of linear fractional transformations with coefficients Gaussian integers. The Picard group is generated by $X'(z) = i/(iz + 1)$, $U'(z) = -1/z$, $Y'(z) = -(z + 1)/z$, and $V'(z) = i/iz$. The Picard group can also be generated by $R'P'$, $S'P'$, $Q'P'Q'P'$ and $Q'R'P'Q'$ where $X' = Q'R'Q'S'$, $U' = S'Q'P'Q'$, $Y' = S'R'$ and $V' = P'S'$ with defining relations $(R'P')^2 = (S'P')^2 = (Q'P'Q'P')^2 = (Q'R'P'Q')^2 = (Q'R'Q'S')^3 = (S'Q'P'Q')^2 = (S'R')^3 = (Q'R'Q'R')^2 = e$. In deriving its subgroups, the Picard group is viewed as a group of symmetries of a three dimensional tiling by a hyperbolic tetrahedron with dihedral angles $\pi/4$, $\pi/4$, $\pi/3$, $\pi/2$, $\pi/2$ and $\pi/2$, having symmetry group generated by $P'$, $Q'$, $R'$ and $S'$.

## 2 An approach in determining subgroups of symmetry groups in hyperbolic space

In our work, the determination of the subgroups of the extended Hecke, Hecke and Picard groups is facilitated by the link of group theory and color symmetry theory. More particularly, the derivation of the subgroups is realized by a correspondence that arises between a subgroup of a group $M$ of symmetries of a given tiling and $M$-transitive colorings of the tiling. This idea is encapsulated in the following result which was established in [4, 13]:

**Theorem 1.** *Let $\mathcal{T}$ be a tiling of Euclidean, spherical or hyperbolic space by copies of a Coxeter polytope D which fill space with no gaps and overlaps. Let H be the group generated by reflections in the facets of D. Let M be a subgroup of H and $\mathcal{O} = \{mD \mid m \in M\}$ be the M-orbit of D.*
(i) *Suppose L is a subgroup of M of index n. Let $\{m_1, m_2, …, m_n\}$ be a complete set of left coset representatives of L in M and $\{c_1, c_2, …, c_n\}$ a set of n colors. The assignment $m_iLD \to c_i$ defines an n-coloring of $\mathcal{O}$ which is M-transitive.*
(ii) *In an M-transitive n-coloring of $\mathcal{O}$, the elements of M which fix a specific color in the colored set $\mathcal{O}$ form a subgroup of M of index n.*

In the above theorem we assume that the tiling $\mathcal{T}$ is the $H$-orbit of $D$, that is $\mathcal{T} = \{hD \mid h \in H\}$; $D$ forms a fundamental region for $H$ and $Stab_H(D) = \{h \in H \mid hD = D\}$ is the trivial group $\{e\}$ [22].

Moreover, for a subgroup $M$ of $H$, $Stab_M(D) = \{e\}$ and $M$ acts transitively on $\mathcal{O} = \{mD \mid m \in M\}$. Consequently, there is a one-to-one correspondence between $M$ and $\mathcal{O}$ given by $m \to mD$, $m \in M$. The action of $M$ on $\mathcal{O}$ is regular, where $m' \in M$ acts on $mD \in \mathcal{O}$ by sending it to its image under $m'$.

The index $n$ subgroups of the extended Hecke group $H$, the Hecke group $K$ and the Picard group $\mathcal{P}$ is derived using $n$-colorings of $\mathcal{O} = MD$ where $M = H$, $K$ and $\mathcal{P}$, respectively. We define an $n$- coloring of the $M$-orbit of $D$ given by $\mathcal{O} = MD = \{mD : m \in M\}$, for a subgroup $M$ of $H$ as follows:

If $C = \{c_1, c_2, …, c_n\}$ is a set of $n$ colors, an onto function $f : \mathcal{O} \to C$ is called an $n$-*coloring* of $\mathcal{O}$. To each $mD \in \mathcal{O}$ is assigned a color in $C$. The coloring determines a partition $\mathcal{P} = \{f^{-1}(c_i) : c_i \in C\}$ where $f^{-1}(c_i)$ is the set of elements of $\mathcal{O}$ assigned color $c_i$. Equivalently, we may think of the coloring as a partition of $\mathcal{O}$.

In this paper, we look at $M$-transitive $n$-colorings of $\mathcal{O}$, that is, $n$-colorings of $\mathcal{O}$ for which the group $M$ has a transitive action on the set $C$ of $n$ colors. For a discussion on how to arrive at such a coloring and the proof of the preceding theorem, the reader may refer to [4].

**Remark:** Note that given a subgroup $L$ of $M$ of index $n$ and a set of $n$ colors $\{c_1, c_2, …, c_n\}$, then there correspond $(n-1)!$ $M$-transitive $n$-colorings of $\mathcal{O}$ with $L$ fixing $c_1$. In an $M$-transitive $n$-coloring of $\mathcal{O}$ with $L$ fixing $c_1$, the set of polytopes $LD$ is assigned $c_1$ and the remaining $n-1$ colors are distributed among the $m_iLD$, $m_i \notin L$.

## 3 Index 2, 3 and 4 subgroups of the extended Hecke and Hecke groups

In this section, the low index subgroups of the extended Hecke and Hecke groups are derived using Theorem 1 with the following underlying assumptions: Consider a hyperbolic triangle $\Delta$ (a two-dimensional Coxeter polytope) with interior angles $\pi/2$, $\pi/p$ and zero. Repeatedly reflecting $\Delta$ in its sides results in a tiling $T$ of hyperbolic 2-space by copies of $\Delta$. Let $P$, $Q$, $R$ denote reflections along the sides of $\Delta$ opposite the angles $\pi/2$, $\pi/p$, 0, respectively. The group $H$ generated by $P$, $Q$, $R$ along the sides of $\Delta$ is called an extended Hecke group and is denoted by *$2p\infty$. The index 2 subgroup of $H$ consisting of orientation-preserving isometries and generated by the reflections $QR$ and $RP$ is the Hecke group $K$. The elements of $H$ are symmetries of $T$ and the triangle $\Delta$ forms a fundamental region for $H$. Note that $Stab_H(\Delta)$, $Stab_K(\Delta)$ is $\{e\}$.

As an illustration, we present in Fig. 1(a) a triangle tiling by a fundamental triangle with interior angles $\pi/2$, $\pi/4$ and zero with the axes of reflections $P$, $Q$, $R$ and centers of rotations $QR$ and $RP$ labeled. In this case, $H = $ *$24\infty$.

Conway's orbifold notation [1] will be used to describe the subgroups. The orbifold notation for crystallographic groups is based on the type of symmetries occurring in the group: * indicates a mirror reflection, $x$ a glide reflection, and a number $n$ indicates a rotation of order $n$. Moreover, if a number $n$ comes after the *, the center of the corresponding rotation lies on mirror lines, so the symmetry there is dihedral of order $2n$.

**Theorem 2.** *The extended Hecke group* $H = *2p\infty$ *has*
(i) 7 *index* 2 *subgroups if p is even* (*Table* 1 *nos.* 1-7).
(ii) 3 *index* 3 *subgroups up to conjugacy in H if p is divisible by both* 2 *and* 3 (*Table* 2 *nos.* 8-10).
(iii) 22 *index* 4 *subgroups up to conjugacy in H if p is divisible by* 2, 3 *and* 4 (*Table* 2 *nos.* 11-32).

**Proof.** To arrive at the index 2, 3 and 4 subgroups of the extended Hecke group *H*, we construct respectively, two, three and four colorings of *T* where all elements of *H* effect color permutations and *H* acts transitively on the set of colors. For such an *n*-coloring of *T*, $n \in \{2, 3, 4\}$, a homomorphism $\pi: H \to S_n$ is defined. The group *H* is generated by *P*, *Q*, and *R*, thus $\pi$ is completely determined when $\pi(P)$, $\pi(Q)$, and $\pi(R)$ are specified. To construct a coloring, we consider a fundamental region for *H* and assign to it color $c_1$. The arguments to obtain 2-, 3- or 4-colorings of *T* where each of the generators *P*, *Q*, *R* permutes the colors, that is, the associated color permutation to a generator is of order 1 or 2, are as follows:

(i) In constructing 2-colorings of *T*, each of *P*, *Q*, *R* may be assigned either the identity permutation (1) or the permutation (12). The possible 2-colorings of *T* are listed in Table 1 line nos. 1-7.

The coloring shown in Fig. 1(a) is a result of the color scheme suggested in Table 1 line no. 1, where each of the reflections *P*, *Q*, *R* interchanges colors $c_1$ and $c_2$. As a consequence of this assignment of colors, the *p*-fold rotation *QR* together with the 2-fold rotation *RP* and the $\infty$-fold rotation *QP* generate the subgroup fixing color $c_1$. The orbifold notation for the group is $p2\infty$.

Consider the coloring shown in Fig. 1(b) corresponding to the color assignment given in Table 1 line no. 2 where *P* interchanges colors $c_1$ and $c_2$ while both *Q* and *R* fix the colors $c_1$ and $c_2$. The reflections *Q*, *R*, *PQP* together generate the subgroup fixing color $c_1$ in the coloring. Note that the products of pairs of these reflections, namely: the *p*-fold rotation *QR*, *p*-fold rotation *PQRP* and $\infty$-fold rotation $(QP)^2$ also fix color $c_1$. The orbifold notation for the given group is $*pp\infty$. Note that the color scheme given in Table 1 line nos. 3 and 4 will result in colorings where the respective subgroups fixing color $c_1$ are groups of similar type.

Shown in Fig. 1(c) is the coloring resulting from the color assignment given in Table 1 line no. 5 where *P* fixes the colors and each of *Q* and *R* interchanges colors $c_1$ and $c_2$. Resulting from this color scheme, the reflection *P* as well as the *p*-fold rotation *QR* generate the subgroup fixing color $c_1$ in the coloring. This group in orbifold notation is $p*\infty$. Groups of similar type result from color assignments given in Table 1 line nos. 6 and 7. Considering the colorings that will yield distinct subgroups there will be 7 subgroups of index 2 in $*2p\infty$.

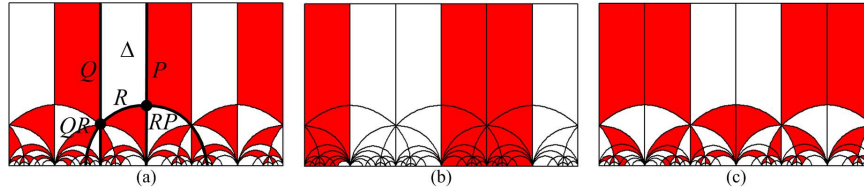

**Fig. 1.** The colorings of the tiling by triangles with angles $\pi/2$, $\pi/4$ and zero resulting from the color assignments given in Table 1: (a) 1 (b) 2 (c) 5. The colors white and red are used to represent, respectively, colors $c_1$ and $c_2$.

(ii) In arriving at 3-colorings of *T*, this time each of *P*, *Q*, *R* may be assigned either the identity permutation or either of the permutations (12), (13) or (23). Now, since we assume *H* acts transitively on the set of 3 colors, we consider the assignments that bring forth a $\pi(H)$ which is a transitive subgroup of $S_3$. Note that there are two transitive subgroups of $S_3$, namely $S_3$ and the cyclic group $Z_3$. It is not possible for $\pi(H)$ to be $Z_3$ since this group cannot be generated by elements of order two. Thus we only consider the case when $\pi(H)$ is $S_3$. In order to generate $S_3$, two of the three permutations (12), (13) or (23) must be associated with two from among the generators *P*, *Q*, and *R*. The remaining third generator is associated either with the identity, with one of the two already assigned earlier or with the third permutation of order 2. The possible 3-colorings that will yield distinct subgroups of index 3 distinct up to conjugacy in *H* are given in Table 1 line nos. 8-10.

The color scheme suggested in Table 1 line no. 8 gives rise to the coloring shown in Fig. 2(a) where the reflection *Q* interchanges $c_1$ and $c_2$ and the reflections *P* and *R* interchange $c_1$ and $c_3$. As a consequence of this assignment of colors, the reflections *QPQ*, *QRQ*, *RQR* together with the 2-fold rotation *RP* generate the subgroup fixing $c_1$. The $(p/3)$-fold rotation $(QR)^3$, the 2-fold rotation *QRPQ* and the $\infty$-fold rotation $(QP)^3$ also fix $c_1$. The orbifold notation for the group is $2*2(p/3)\infty$.

The coloring given in Fig. 2(b) corresponding to the color assignment given in Table 1 line no. 9 has *P* fixing all the colors, *Q* interchanging $c_1$ and $c_2$ and *R* interchanging $c_1$ and $c_3$. The reflections *QRQ*, *RQR*, *P* and *QPQ* generate the subgroup fixing $c_1$ of the given coloring; where the region bounded by the axes of these reflections serves as a fundamental region. Note that the products of pairs of these reflections, namely: the $(p/3)$-fold rotation $(QR)^3$, the $\infty$-fold rotation *RPQR*, the $\infty$-fold rotation $(PQ)^2$ and the 2-fold rotation *QRPQ* also fix $c_1$. In orbifold notation, the group is $*2(p/3)\infty\infty$. Note that the color assignment given in Table 1 line no. 10 will result in a coloring where the respective subgroup fixing $c_1$ is a group of similar type.

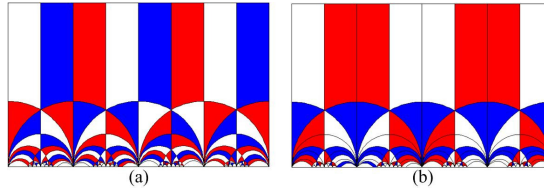

**Fig. 2.** The colorings of the tiling by triangles with angles $\pi/2$, $\pi/3$ and zero resulting from the color assignments given in Table 1: (a) 8 (b) 9. The colors white, red and blue are used to represent, respectively, colors $c_1$, $c_2$ and $c_3$.

(iii) In constructing 4-colorings of $T$ that will yield index 4 subgroups of $H$, each of $P$, $Q$, $R$ is assigned either the identity permutation; a 2-cycle or a product of two disjoint 2-cycles: (12), (13), (14), (23), (24), (34), (12)(34), (13)(24) or (14)(23). We only look at the situations that arise when $\pi(H)$ is a transitive subgroup of $S_4$. In particular, we consider the cases when $\pi(H)$ is either $S_4$, $D_4$ or $V = \{e, (12)(34), (13)(24), (14)(23)\}$, a Klein 4 group. The possible $H$-transitive 4-colorings of $T$ that can be constructed when $\pi(H) \cong V$, $D_4$ and $S_4$ are given in Table 1 line nos. 11-32. These colorings will give rise to subgroups of index 4 distinct up to conjugacy in $H$.

Consider the coloring given in Fig. 3(a) corresponding to the color assignment given in Table 1 line no. 11 where the reflections $P$, $Q$, $R$ respectively interchange $c_2$ and $c_4$; $c_1$ and $c_2$; $c_1$ and $c_3$. Consequently, the reflections $P$, $QPQPQ$, $QRQ$, $RQR$ generate the subgroup fixing $c_1$ of the given coloring; where the region bounded by the axes of reflections serves as a fundamental region. Note that the products of pairs of these reflections, namely: the $(p/3)$-fold rotation $(QR)^3$, the $\infty$-fold rotation $(PQ)^3$, the $p$-fold rotation $QPQRPQ$ and the $\infty$-fold rotation $RQPR$ also fix $c_1$. In orbifold notation, the group is $*p(p/3)\infty\infty$. Similar groups will fix $c_1$ of the colorings resulting from the color assignments from Table 1 line nos. 14-17.

For the given coloring appearing in Fig. 3(b), corresponding to the color assignment in Table 1 line no. 12, the reflection $R$ interchanges $c_1$ and $c_2$; $c_3$ and $c_4$; whereas $P$ interchanges the colors $c_1$ and $c_2$; and $Q$ interchanges the colors $c_1$ and $c_3$. As a result, the reflections $RQR$, $QRQRQ$, $QPQ$ together with the 2-fold rotation $RP$ generate the subgroup fixing $c_1$ of the given coloring. The $(p/4)$-fold rotation $(QR)^4$, the $\infty$-fold rotation $QRPQRQ$ and the $\infty$-fold rotation $(PQ)^3$ also fix $c_1$. The orbifold notation for the given group is $2*(p/4)\infty\infty$. Note that the color assignments given in Table 1 line nos. 13, 18, and 27-29 will give rise to groups of similar type.

Now, other types of subgroups fixing color $c_1$ can be derived from the remaining color assignments given in Table 1. The coloring shown in Fig. 3(c) results from the color scheme presented in Table 1 line no. 19 where $c_1$ and $c_3$; $c_2$ and $c_4$ are at the same time interchanged by the reflection $Q$ while the reflections $R$ and $P$ interchange $c_1$ and $c_2$; $c_3$ and $c_4$, respectively. As a consequence of this assignment of colors, the reflections $P$ and $QRQ$ fix $c_1$ in the given coloring. The translation $RQPQ$ also fixes $c_1$. Together they generate the subgroup fixing $c_1$ of the corresponding coloring. In addition, the $\infty$-fold rotation $(PQ)^4$ and the $(p/4)$-fold rotation $(QR)^4$ also fix $c_1$ of the given coloring. In orbifold notation, this subgroup is $*\infty*(p/4)$.

Given in Fig. 3(d) is a coloring arising from the color scheme in Table 1 line no. 20, where the reflection $R$ interchanges $c_1$ and $c_2$; $c_3$ and $c_4$; whereas $P$ interchanges the colors $c_1$ and $c_2$; and $Q$ interchanges the colors $c_1$ and $c_3$; $c_2$ and $c_4$. The reflections $QPQ$ and $QRPRQ$ together with the 2-fold rotation $RP$ and $(p/2)$-fold rotation $(QR)^2$ generate the subgroup fixing $c_1$ of the given coloring. Also fixing $c_1$ is the $\infty$-fold rotation $(QP)^4$. The orbifold notation for the given group is $(p/2)2*\infty$. Note that the color assignments given in Table 1 line nos. 21-24 will give rise to groups of similar type.

Fig. 3(e), on the other hand, shows a coloring of the tiling obtained from the color assignment given in Table 1 line no. 25 where $c_1$ and $c_4$; $c_2$ and $c_3$ are interchanged by the reflection $P$ while the reflection $Q$ interchanges $c_1$ and $c_2$; $R$ interchanges $c_1$ and $c_3$; $c_2$ and $c_4$. As a consequence of this assignment of colors, the reflections $RQR$ and $PQP$ together with the glide reflection $QPR$ generate the subgroup fixing $c_1$ of the given coloring. The $(p/4)$-fold rotation $(QR)^4$ and $\infty$-fold rotation $(PQ)^4$ also fix $c_1$. In orbifold notation, this group is $*(p/4)\infty x$.

The color scheme suggested in line no. 26, Table 1, gives rise to the coloring in Fig. 3(f), where the subgroup fixing $c_1$ is generated by three glide reflections, namely: $PQR$, $PRQ$, $QPR$. Note that products of pairs of glide reflections result in the following rotations fixing $c_1$: the $(p/2)$-fold rotation $(QR)^2$ and the $\infty$-fold rotation $(PQ)^2$. The orbifold notation for this subgroup is $(p/2)\infty x$.

Finally, consider the coloring in Fig. 3(g) resulting from the color assignment given in Table 1 line no. 30. In this coloring, three rotations generate the subgroup fixing $c_1$. These are the two $p$-fold rotations $QR$, $PQRP$ and the $\infty$-fold rotation $(PQ)^2$. This subgroup in orbifold notation is $pp\infty$. Groups of similar type arise from the color assignment given in Table 1 line nos. 31-32. This completes the derivation of the 22 index 4 subgroups of $H$ up to conjugacy in $H$. ∎

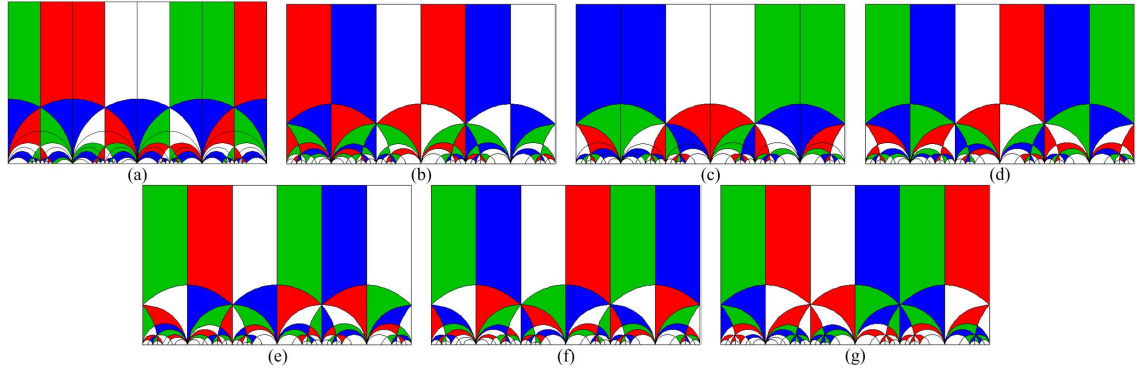

**Fig. 3.** The colorings of the tiling by triangles with angles (a) π/2, π/3 and zero; (b)-(g) π/2, π/4 and zero resulting from the color assignments given in Table 1: (a) 11; (b) 12; (c) 19; (d) 20; (e) 25; (f) 26; (g) 30. The colors white, red, blue and green are used to represent, respectively, colors $c_1$, $c_2$, $c_3$ and $c_4$.

Table 1. Colorings that will give rise to index 2, 3, 4 subgroups of $H = *2p\infty$.

| no | index in $H$ | $P$ | $Q$ | $R$ | generators | a fundamental region | orbifold notation |
|---|---|---|---|---|---|---|---|
| 1 | 2 | (12) | (12) | (12) | QR, RP, QP | $\Delta \cup Q(\Delta)$ | $2p\infty$ |
| 2 | 2 | (12) | (1) | (1) | R, Q, PQP | $\Delta \cup P(\Delta)$ | $*pp\infty$ |
| 3 | 2 | (1) | (12) | (1) | P, R, QRQ, QPQ | $\Delta \cup Q(\Delta)$ | $*(p/2)\infty\infty$ |
| 4 | 2 | (1) | (1) | (12) | RQR, P, Q | $\Delta \cup R(\Delta)$ | $*22(p/2)\infty$ |
| 5 | 2 | (1) | (12) | (12) | RQ, P | $\Delta \cup R(\Delta)$ | $p*\infty$ |
| 6 | 2 | (12) | (1) | (12) | RP, Q, PQP | $\Delta \cup P(\Delta)$ | $\infty*(p/2)$ |
| 7 | 2 | (12) | (12) | (1) | QP, R, QRQ | $\Delta \cup Q(\Delta)$ | $2*(p/2)\infty$ |
| 8 | 3 | (13) | (12) | (13) | RP, QPQ, QRQ, RQR | $\Delta \cup Q(\Delta) \cup R(\Delta)$ | $2*2(p/3)\infty$ |
| 9 | 3 | (1) | (12) | (13) | QRQ, RQR, P, QPQ | $\Delta \cup Q(\Delta) \cup R(\Delta)$ | $*2(p/3)\infty\infty$ |
| 10 | 3 | (12) | (13) | (1) | PQP, QPQ, QRQ, R | $\Delta \cup P(\Delta) \cup Q(\Delta)$ | $*2(p/2)p\infty$ |
| 11 | 4 | (24) | (12) | (13) | P, QPQPQ, QRQ, RQR | $\Delta \cup R(\Delta) \cup Q(\Delta) \cup QP(\Delta)$ | $*p(p/3)\infty\infty$ |
| 12 | 4 | (12) | (13) | (12)(34) | RP, RQR, QRQRQ, QPQ | $\Delta \cup Q(\Delta) \cup R(\Delta) \cup QR(\Delta)$ | $2*(p/4)\infty\infty$ |
| 13 | 4 | (13)(24) | (12) | (13) | RP, PQP, QPQPQ, QRQ | $\Delta \cup Q(\Delta) \cup P(\Delta) \cup QP(\Delta)$ | $2*p(p/3)\infty$ |
| 14 | 4 | (1) | (12) | (13)(24) | QRQRQ, RQR, P, QPQ | $\Delta \cup Q(\Delta) \cup R(\Delta) \cup QR(\Delta)$ | $*(p/4)\infty\infty\infty\infty$ |
| 15 | 4 | (12) | (13)(24) | (1) | PQPQP, QPQ, QRQ, R, PQRQP | $\Delta \cup P(\Delta) \cup Q(\Delta) \cup PQ(\Delta)$ | $**\infty 2(p/2)(p/2)2$ |
| 16 | 4 | (1) | (13)(24) | (12) | RQRQR, QRQ, QPQ, P, RQPQR | $\Delta \cup R(\Delta) \cup Q(\Delta) \cup RQ(\Delta)$ | $*2(p/4)2\infty\infty$ |
| 17 | 4 | (13)(24) | (12) | (1) | QPQPQ, PQP, R, QRQ | $\Delta \cup Q(\Delta) \cup P(\Delta) \cup QP(\Delta)$ | $*p(p/2)p\infty$ |
| 18 | 4 | (12) | (13)(24) | (12) | RP, QRQ, QPQ, PQPQP, PQRQP | $\Delta \cup P(\Delta) \cup Q(\Delta) \cup PQ(\Delta)$ | $2*2(p/4)2\infty$ |
| 19 | 4 | (34) | (13)(24) | (12) | P, QRQ, RQPQ | $\Delta \cup Q(\Delta) \cup R(\Delta) \cup QP(\Delta)$ | $*\infty\infty(p/4)$ |
| 20 | 4 | (12) | (13)(24) | (12)(34) | QRQR, RP, QPQ, QRPRQ | $\Delta \cup Q(\Delta) \cup R(\Delta) \cup QR(\Delta)$ | $(p/2)2*\infty$ |
| 21 | 4 | (13)(24) | (12) | (13)(24) | PQ, RQR, RPQPR | $\Delta \cup R(\Delta) \cup P(\Delta) \cup RP(\Delta)$ | $\infty*(p/4)\infty$ |
| 22 | 4 | (13)(24) | (12) | (12)(34) | QR, PQP, PRQRP | $\Delta \cup P(\Delta) \cup R(\Delta) \cup PR(\Delta)$ | $p*(p/2)\infty$ |
| 23 | 4 | (12)(34) | (13)(24) | (12) | PQPQ, RP, QRQ | $\Delta \cup Q(\Delta) \cup P(\Delta) \cup QP(\Delta)$ | $2\infty*(p/4)$ |
| 24 | 4 | (13)(24) | (12) | (13)(24) | RP, QRPQ, QPQPQ, PQP | $\Delta \cup Q(\Delta) \cup P(\Delta) \cup QP(\Delta)$ | $22*(p/4)\infty$ |
| 25 | 4 | (14)(23) | (12) | (13)(24) | RQR, PQP, QRP | $\Delta \cup Q(\Delta) \cup R(\Delta) \cup PR(\Delta)$ | $*(p/4)\infty x$ |
| 26 | 4 | (12)(34) | (13)(24) | (14)(23) | PQR, PRQ, QPR | $\Delta \cup P(\Delta) \cup Q(\Delta) \cup R(\Delta)$ | $(p/2)\infty x$ |
| 27 | 4 | (1) | (12)(34) | (13)(24) | QRQR, P, QPQ | $\Delta \cup Q(\Delta) \cup R(\Delta) \cup QR(\Delta)$ | $(p/2)*\infty\infty$ |
| 28 | 4 | (13)(24) | (1) | (12)(34) | RQR, Q, PQP, PRQRP | $\Delta \cup P(\Delta) \cup R(\Delta) \cup PR(\Delta)$ | $*(p/2)\infty(p/2)\infty$ |
| 29 | 4 | (12)(34) | (13)(24) | (1) | PQPQ, QRQ, R, PQRQP | $\Delta \cup P(\Delta) \cup Q(\Delta) \cup PQ(\Delta)$ | $\infty*(p/2)(p/2)$ |
| 30 | 4 | (13)(24) | (12)(34) | (12)(34) | QR, PQRP, PQPQ | $\Delta \cup P(\Delta) \cup Q(\Delta) \cup PQ(\Delta)$ | $pp\infty$ |
| 31 | 4 | (13)(24) | (13)(24) | (12)(34) | RP, QRPQ, PQPQ | $\Delta \cup Q(\Delta) \cup P(\Delta) \cup QP(\Delta)$ | $2(p/2)2\infty$ |
| 32 | 4 | (12)(34) | (12)(34) | (13)(24) | PQ, RPQR | $\Delta \cup R(\Delta) \cup P(\Delta) \cup RP(\Delta)$ | $(p/2)\infty\infty$ |

**Theorem 3.** *The Hecke group $K = 2p\infty$ has*
(i) 3 *index* 2 *subgroups if p is even* (*Table* 2 *nos.* 1-3).
(ii) 3 *index* 3 *subgroups up to conjugacy in K if p is divisible by both* 2 *and* 3 (*Table* 2 *nos.* 4-6).
(iii) 10 *index* 4 *subgroups up to conjugacy in K if p is divisible by* 2*,* 3 *and* 4 (*Table* 2 *nos.* 7-16).

**Proof.** In determining the subgroups of the Hecke group $K = 2p\infty$, we consider the $K$-orbit of $\Delta$, namely $K\Delta = \{k\Delta : k \in K\}$ and construct $K$-transitive colorings of $K\Delta$. In determining respectively, the index 2, 3 and 4 subgroups of $K$, we construct two, three and four colorings of $K\Delta$. For an $n$-coloring of $K\Delta$, $n \in \{2, 3, 4\}$, a homomorphism $\pi': K \to S_n$ is defined. The group $K$ is generated by $QR$ and $RP$, thus $\pi'$ is completely determined when $\pi'(QR)$ and $\pi'(RP)$ are specified. To obtain 2-, 3- or 4-colorings of $K\Delta$ where each of the generators $QR$ and $RP$ permutes the colors, the orders of the color permutations associated to $QR$ and $RP$ should be divisors respectively of the orders of $QR$ and $RP$. Construction of a table similar to Table 1, this time pertaining to assignments of permutations to *QR* and *RP*, will yield transitive two, three and four colorings of *K*Δ, listed in Table 2.

(i) To obtain 2-colorings of *K*Δ that will yield index 2 subgroups of *K*, we assign to *QR* and *RP* either the identity permutation (1) or the 2-cycle (12). The consequent color assignments will give rise to three index 2 subgroups of *K* which are precisely the index 4 subgroups of the extended Hecke group *2*p*∞ (Table 1 nos. 30-32) generated by rotations. These groups in orbifold notation are: $pp\infty$, $2(p/2)2\infty$ and $(p/2)\infty\infty$.

(ii) In arriving at 3-colorings of *K*Δ that will give rise to index 3 subgroups of *K*, we assign to *QR* and *RP* color permutations that will generate a transitive subgroup of $S_3$. We will only consider the permutation assignments that will result in $\pi'(K)$ isomorphic to $Z_3$ or $S_3$. The color assignments that give rise to index 3 subgroups of *K* are listed in Table 2 line nos. 4-6. Associated with the color scheme in line no. 4 (see Fig. 4(a) for the corresponding 3-coloring of *K*Δ) is the index 3 subgroup fixing color $c_1$ generated by three 2-fold rotations and one ∞-fold rotation, namely: *RP, QRPQ, PQRPQP* and $(QP)^3$ respectively. Also fixing color $c_1$ is the $(p/3)$-fold rotation $QP(QR)^3PQ$. The group in Conway notation is $(p/3)222\infty$.

Resulting from the color scheme suggested in line no. 5 is the subgroup fixing color $c_1$ generated by one 2-fold rotation and two ∞-fold rotations. These rotations are *RP, QRQP,* $(QP)^3$, respectively. Note the $(p/2)$-fold rotation $(QR)^2$ also fixes color $c_1$. The four rotations generate the index 3 subgroup of *K* with Conway notation $2(p/2)\infty\infty$. Corresponding to this subgroup is the 3-coloring of *K*Δ shown in Fig. 4(b). The subgroup that fixes $c_1$ pertaining to the color assignment given in Table 2 line no. 6 is also of same type.

(iii) In determining the 4-colorings of *K*Δ that will correspond to index 4 subgroups of *K*, we assign to *QR* and *RP* color permutations that will yield any transitive subgroup of $S_4$. We will consider the permutation assignments that bring forth a $\pi'(K)$ isomorphic to $V = \{e, (12)(34), (13)(24), (14)(23)\}$, $Z_4$, $D_4$, $A_4$, or $S_4$. The possible *K*-transitive 4-colorings of *K*Δ that can be constructed when $\pi'(K)$ is isomorphic to these groups are given in Table 2 line nos. 7-16. These colorings will give rise to subgroups of index 4 distinct up to conjugacy in *K* which we characterize in terms of the Conway orbifold notation as follows:

The color scheme in line no. 7 results in the 4-coloring of *K*Δ (Fig. 4(c)) where the ∞-fold rotation $(QP)^2$ together with two $(p/2)$-fold rotations $(RQ)^2$, $P(RQ)^2P$ fix $c_1$ and generate the corresponding index 4 subgroup. In Conway notation, this group is $(p/2)(p/2)\infty$.

The 4-coloring of *K*Δ (Fig. 4(d)) resulting from the color assignment given in line no. 8 corresponds to the subgroup fixing $c_1$ generated by three 2-fold rotations, namely, *RP, QRPQ, PQRPQP* and two ∞-fold rotations $(QP)^4$, *QPQRQPQP*. Also fixing $c_1$ in the coloring is the $(p/4)$-fold rotation $P(RQ)^4P$. The orbifold notation for the corresponding index 4 subgroup is $222(p/4)\infty\infty$.

Moreover, the 4-coloring of *K*Δ (Fig. 4(e)) obtained from the color scheme in line no. 9 gives rise to the subgroup fixing $c_1$ generated by three ∞-fold rotations, namely: *QRQP, RQPQ, PRQPQP*. Note that the $(p/4)$-fold rotation $P(RQ)^4P$ also fixes $c_1$. The orbifold notation for the associated index 4 subgroup is $(p/4)\infty\infty\infty\infty$.

The color scheme suggested in line no. 10 results in the 4-coloring of *K*Δ (Fig. 4(f)) where the color $c_1$ is fixed by *RP* which is 2-fold; *QRQP,* $(QP)^4$, *QPQRQPQP* which are all ∞-fold. In addition, $(QR)^2$, a $(p/2)$-fold, also fixes $c_1$. The index 4 subgroup in Conway notation is $2(p/2)\infty\infty\infty$.

The subgroup generated by the 2-fold *RP* and ∞-fold rotations $(QP)^2$, *QRQPQRQP, QRQRQRQP* corresponds to the coloring given in line no. 11 (Fig. 4(g)). In the same coloring, $c_1$ is also fixed by the $(p/4)$-fold rotation $(QR)^4$. This group in orbifold notation is $2(p/4)\infty\infty\infty$. Subgroups of similar type arise from the colorings listed in Table 2 line nos. 12-13.

Now, the color scheme given in line no. 14 gives rise to the coloring of *K*Δ (Fig. 4(h)) where $c_1$ is fixed by two 2-fold rotations *RQRPQR, PRQRPQRP* and two ∞-fold rotations *RQPQPQRP, QP*. In orbifold notation, the group is $22\infty\infty$. The color assignment in line no. 15 corresponds to the subgroup fixing $c_1$ in the coloring (Fig. 4(g)) generated by two 2-fold rotations *RP, QR* and a *p*-fold rotation *PQPRQPQP* and ∞-fold *QPQRQP*. The $(p/3)$-fold $(RQ)^3$ also fixes $c_1$. The group in orbifold notation is $22(p/3)p\infty$. A similar type of group arises from the color assignment given in the last line in Table 2. ∎

**Corollary 4.** *The modular group K' = 23∞ has 1 index 2 subgroup, 2 index 3 subgroups and 2 index 4 subgroups.*

**Proof.** This follows immediately from Theorem 3 and the fact that *QR* is an element of order 3. Refer to Table 2 for the colorings that give rise to the low index subgroups: index 2 (Table 2 no. 1); index 3 (Table 2 nos. 4 and 6) and index 4 (Table 2 nos. 14 and 15). ∎

Table 2. Colorings that will give rise to index 2, 3, 4 subgroups of $K = 2p\infty$. Here, $\Delta' = \Delta \cup Q(\Delta)$.

| no | index in $K$ | QR | RP | QP | generators | a fundamental region | orbifold notation |
|---|---|---|---|---|---|---|---|
| 1 | 2 | (1) | (12) | (12) | QR, PQRP, PQPQ | $\Delta' \cup PQ(\Delta')$ | $pp\infty$ |
| 2 | 2 | (12) | (1) | (12) | RP, QRPQ, PQPQ | $\Delta' \cup QP(\Delta')$ | $22(p/2)\infty$ |
| 3 | 2 | (12) | (12) | (1) | PQ, RQPR | $\Delta' \cup RQ(\Delta')$ | $(p/2)\infty\infty$ |
| 4 | 3 | (123) | (1) | (123) | RP, QRPQ, PQRQPQ, QPQPQP | $\Delta' \cup PQ(\Delta') \cup QP(\Delta')$ | $(p/3)222\infty$ |
| 5 | 3 | (12) | (23) | (123) | RP, QRQP, QPQPQP | $\Delta' \cup QP(\Delta') \cup PQ(\Delta')$ | $(p/2)2\infty\infty\infty$ |
| 6 | 3 | (123) | (13) | (23) | QP, RQPRP, RQRQRP | $\Delta' \cup RQ(\Delta') \cup QR(\Delta')$ | $2(p/3)\infty\infty\infty\infty\infty$ |
| 7 | 4 | (12)(34) | (13)(24) | (14)(23) | QPQP, RQRQ, PRQRQP | $\Delta' \cup PQ(\Delta') \cup RQ(\Delta') \cup RP(\Delta')$ | $(p/2)(p/2)\infty$ |
| 8 | 4 | (1234) | (1) | (1234) | RP, QRPQ, PQRQP, QPQPQPQP, QPQRQPQP | $\Delta' \cup PQ(\Delta') \cup QP(\Delta') \cup (PQ)^2(\Delta')$ | $222(p/4)\infty\infty\infty\infty$ |
| 9 | 4 | (1234) | (13)(24) | (1432) | QRQP, RQPQ, PRQPQP | $\Delta' \cup PQ(\Delta') \cup RQ(\Delta') \cup RP(\Delta')$ | $(p/4)\infty\infty\infty\infty\infty$ |
| 10 | 4 | (12)(34) | (23) | (1243) | RP, QRQP, QPQPQPQP, QPQRQPQP | $\Delta' \cup PQ(\Delta') \cup QP(\Delta') \cup (PQ)^2(\Delta')$ | $2(p/2)\infty\infty\infty\infty\infty$ |
| 11 | 4 | (1234) | (24) | (12)(34) | RP, QPQP, QRQPQRQP, QRQRQRQR | $\Delta' \cup QR(\Delta') \cup RQ(\Delta') \cup (RQ)^2(\Delta')$ | $2(p/4)\infty\infty\infty\infty\infty$ |
| 12 | 4 | (1243) | (13)(24) | (23) | QP, RQPRP, RQRQRP | $\Delta' \cup QR(\Delta') \cup RQ(\Delta') \cup (RQ)^2(\Delta')$ | $(p/4)\infty\infty\infty\infty\infty\infty$ |
| 13 | 4 | (23) | (12)(34) | (1342) | RQ, PRQRP, PQPQPQP | $\Delta' \cup PQ(\Delta') \cup QP(\Delta') \cup (PQ)^2(\Delta')$ | $(p/2)p\infty\infty\infty$ |
| 14 | 4 | (123) | (13)(24) | (243) | QP, RQRPRP, RQPQPQRPRP | $\Delta' \cup PQ(\Delta') \cup QP(\Delta') \cup (PQ)^2(\Delta')$ | $(p/3)\infty\infty\infty\infty\infty\infty\infty$ |
| 15 | 4 | (123) | (24) | (1243) | RP, QRPQ, QPQPQP, PQPRQPQP | $\Delta' \cup PQ(\Delta') \cup QP(\Delta') \cup (PQ)^2(\Delta')$ | $22(p/3)p\infty$ |
| 16 | 4 | (1234) | (14) | (234) | QP, RQRPQR, PRQRPQRP, RQPQPQRP | $\Delta' \cup PQ(\Delta') \cup QP(\Delta') \cup (PQ)^2(\Delta')$ | $22(p/4)\infty\infty\infty$ |

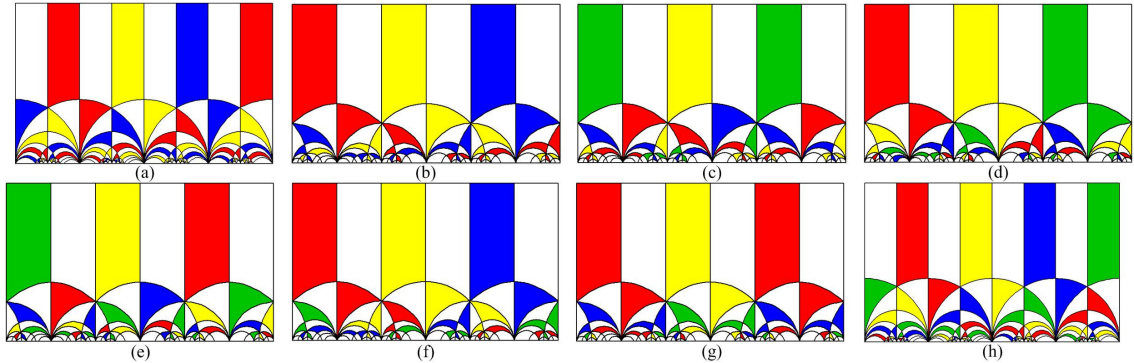

**Fig. 4.** The colorings of the $K$-orbit of a triangle with angles (a) and (h) $\pi/2$, $\pi/3$ and zero; (b)-(g) $\pi/2$, $\pi/4$ and zero. These result from the color assignments given in Table 2: (a) 4; (b) 5; (c) 7; (d) 8; (e) 9; (f) 10; (g) 11; (h) 14. The colors yellow, red, blue and green are used to represent, respectively, colors $c_1$, $c_2$, $c_3$ and $c_4$.

## 4 Index 2, 3 and 4 subgroups of the Picard group

In deriving subgroups of the Picard group $\mathcal{P}$, our starting point would be to consider the hyperbolic tetrahedron t (a three dimensional Coxeter polytope) with dihedral angles $\pi/4$, $\pi/4$, $\pi/3$, $\pi/2$, $\pi/2$, $\pi/2$. Repeatedly reflecting t in its faces results in a three-dimensional hyperbolic tiling $\mathcal{T}$ by copies of t. The group $G$ generated by the reflections $P'$, $Q'$, $R'$, $S'$ along the respective faces of t is a tetrahedron group with one of its index 4 subgroups the Picard group $\mathcal{P}$. $G$ is the symmetry group of $\mathcal{T}$ and t forms a fundamental region for $G$. $\mathcal{P}$ is transitive on the $\mathcal{P}$-orbit of t, given by $\mathcal{P}$t and $Stab_{\mathcal{P}}(\text{t})$ is $\{e\}$. The following result explains the derivation of the low index subgroups of $\mathcal{P}$ given this setting and using Theorem 1:

**Theorem 5.** *The Picard group $\mathcal{P}$ has*
(i) 3 *index* 2 *subgroups* (*Table* 3 *nos*. 1-3).
(ii) 2 *index* 3 *subgroup up to conjugacy in $\mathcal{P}$* (*Table* 3 *nos*. 4-5).
(iii) 7 *index 4 subgroups up to conjugacy in $\mathcal{P}$* (*Table* 3 *nos*. 6-12).

**Proof.** In deriving the index 2, 3, 4 subgroups of $\mathcal{P}$, we construct $\mathcal{P}$-transitive 2-, 3-, 4-colorings of $\mathcal{P}$t. For each coloring, we define a homomorphism $\pi''$: $\mathcal{P} \to S_n$, $n \in \{2, 3, 4\}$. This homomorphism is completely determined when $\pi''(R'P')$, $\pi''(S'P')$, $\pi''(Q'P'Q'P')$, and $\pi''(Q'R'P'Q')$, the images of the generators of $\mathcal{P}$ under the homomorphism, are specified. Following the arguments made in Theorems 1 and 2, we will only look at the cases where $\pi''(\mathcal{P}\text{t})$ is a transitive subgroup of $S_n$. Moreover, since the rotations $Q'R'Q'S'$, $S'R'$ are order three elements; $S'Q'P'Q'$ and $Q'R'Q'R'$ are order two elements, we only consider the permutation assignments to $R'P'$, $S'P'$, $Q'P'Q'P'$, $Q'R'P'Q'$ resulting in permutations corresponding to $Q'R'Q'S'$ and $S'R'$ that are 3-cycles and to $S'Q'P'Q'$ and $Q'R'Q'R'$ that are 2-cycles. This gives us three, two, and seven possible color assignments which correspond respectively to the index 2, 3, 4 subgroups distinct up to conjugacy in $\mathcal{P}$. ∎

Table 3. Colorings that will give rise to index 2, 3, 4 subgroups of the Picard group 𝒫.

| no | index in 𝒫 | R'P' | S'P' | Q'P'Q'P' | Q'R'P'Q' | Q'R'Q'S' | S'Q'P'Q' | S'R' | Q'R'Q'R' | generators |
|---|---|---|---|---|---|---|---|---|---|---|
| 1 | 2 | (1) | (1) | (12) | (1) | (1) | (12) | (1) | (1) | P'S', P'Q'P'R'Q'P', P'Q'R'Q'P'Q'P' |
| 2 | 2 | (12) | (12) | (1) | (12) | (1) | (12) | (1) | (1) | R'S', P'Q'P'Q', Q'R'P'Q'S'P' |
| 3 | 2 | (12) | (12) | (12) | (12) | (1) | (1) | (1) | (1) | S'R', Q'P'Q'R', Q'R'P'Q'R'P' |
| 4 | 3 | (23) | (12) | (1) | (23) | (132) | (12) | (123) | (1) | R'P', Q'P'Q'P', Q'R'P'Q', S'R'Q'R'P'Q'S'P', P'S'Q'R'P'Q'S'P'Q'R'P'Q'S'P' |
| 5 | 3 | (23) | (12) | (12) | (23) | (132) | (1) | (123) | (1) | R'P', Q'R'P'Q', Q'P'S'Q', S'R'S'P'R'S', S'R'Q'P'Q'P'R'S', P'S'Q'R'P'Q'R'S', S'R'Q'R'P'Q'S'P', P'S'Q'R'Q'P'Q'R'P'Q'S'P', P'S'Q'R'P'Q'S'P'Q'R'P'Q'S'P' |
| 6 | 4 | (12)(34) | (12)(34) | (13)(24) | (12)(34) | (1) | (14)(23) | (1) | (1) | R'S', Q'R'Q'S', Q'R'P'Q'S'P', S'Q'P'R'Q'P', S'Q'R'Q'R'Q'P' |
| 7 | 4 | (34) | (23) | (1) | (12) | (123) | (23) | (234) | (12)(34) | R'P', S'P', Q'P'Q'P', Q'R'Q'P'Q'R'P'Q', Q'R'Q'S'P'Q'R'P'Q', Q'R'P'Q'R'P'Q'R'P'Q', Q'R'P'Q'S'R'Q'P'Q'R'S'P'R'P'Q', Q'R'P'Q'S'R'Q'R'Q'R'S'Q'R'P'Q', Q'R'P'Q'S'R'Q'R'P'Q'P'S'Q'R'P'Q' |
| 8 | 4 | (34) | (23) | (23) | (12) | (123) | (1) | (234) | (12)(34) | R'P', S'P', Q'P'Q'P', Q'R'Q'P'Q'R'Q', Q'R'P'Q'R'P'Q'R'P', Q'R'Q'S'R'Q'P'R'S'Q'R'P'Q' |
| 9 | 4 | (34) | (23) | (14) | (34) | (243) | (14)(23) | (234) | (1) | R'P', S'P', Q'R'Q'P', P'Q'P'Q'S'Q'P'Q', P'Q'R'Q'P'S'R'S'Q'R'Q'P', P'Q'R'P'S'Q'P'Q'S'Q'R'Q'P', P'Q'R'Q'P'S'Q'P'R'Q'S'P'Q'R'Q'P', P'Q'P'R'Q'P'R'Q'P', P'Q'P'R'P'Q'R'Q'P' |
| 10 | 4 | (34) | (23) | (14) | (12) | (123) | (14)(23) | (234) | (12)(34) | R'P', S'P', Q'P'R'Q'S'R'Q'P'Q'P', Q'P'Q'R'Q'P'Q'R'Q'P'Q'P', Q'P'Q'R'Q'P'R'Q'R'Q'P'Q', P'Q'P'Q'S'Q'P'Q', P'Q'R'P'Q'P', Q'P'R'Q'R'P'Q'P'R'Q', P'Q'P'R'S'Q'P'R'Q', Q'P'R'Q'P'Q'R'Q' |
| 11 | 4 | (34) | (23) | (14)(23) | (34) | (243) | (14) | (234) | (1) | R'P', S'P', Q'R'Q'P', Q'P'Q'R'Q'R'Q'P', P'Q'R'Q'P'S'R'Q'R'Q'P', P'Q'R'P'S'R'Q'R'Q'R'Q'P', P'Q'R'P'S'R'Q'P', P'Q'R'P'S'R'Q'P', P'Q'R'Q'P'S'Q'R'P'Q'S'P'Q'R'Q'P' |
| 12 | 4 | (34) | (23) | (14)(23) | (12) | (123) | (14) | (234) | (12)(34) | R'P', S'P', P'Q'P'R'Q'P', P'Q'P'Q'S'Q'P'Q', Q'R'Q'R'Q'P'Q'P', Q'P'Q'R'Q'R'P'Q'R'Q'P'Q', Q'P'Q'R'S'Q'P'Q'R'Q'P'Q', P'Q'R'P'Q'R'Q'P'Q'R'Q'P'Q' |

## 5 Outlook

As a continuation of this work, it would be interesting to identify the embeddings of the Hecke groups in the Picard group [6]. There is still a lot to be investigated on the subgroup structure of the extended Hecke and Hecke groups. In [12, 19], the commutator subgroups, even subgroups, and power subgroups of the extended Hecke group have been studied. The conjugacy classes of torsion elements of the extended Hecke group have been investigated in [17], and this information has been used to find some normal subgroups and Fuchsian subgroups. It would be worthwhile to verify and generalize these results using the methods found in this paper. Moreover, it would be useful to obtain similar results for the Picard group [16, 20].

The next step would be to study the subgroup structure of higher dimensional analogues of the modular and Picard groups. In [14], Ahlfor's description of hyperbolic isometries via Clifford algebras $\mathscr{C}_n$[1, 2] gives rise to a natural generalization of the modular and Picard groups in higher dimensions, such as $PSL(2, Z[i, j])$. Studying these analogues of the modular and Picard groups using tools from color symmetry theory would be noteworthy. Future work would also include the study of quaternionic modular groups, which occur as subgroups of groups generated by reflections with axes along the facets of four- and five-dimensional hyperbolic Coxeter polytopes [11].

**Acknowledgements.** Ma. Louise De Las Peñas would like to acknowledge funding support of the Ateneo de Manila University Loyola Schools through the Scholarly Work Faculty Grant. Ma. Carlota Decena is grateful for the support from the Research Center for Natural Sciences, University of Santo Tomas.